\documentclass{amsart}
\usepackage{amssymb}
\usepackage{bbm}
\usepackage{amsfonts}
\usepackage{mathrsfs}

\newtheorem{theorem}{Theorem}[section]
\newtheorem{lemma}[theorem]{Lemma}

\theoremstyle{definition}
\newtheorem{definition}[theorem]{Definition}

\newtheorem{cor}[theorem]{Corollary}
\theoremstyle{remark}
\newtheorem{remark}[theorem]{Remark}

\numberwithin{equation}{section}

\def\and{\cap}

\def\bref#1{(\ref{#1})}
\def\proof{{\noindent\em Proof:} }

\DeclareFontFamily{U}{fsy}{} \DeclareFontShape{U}{fsy}{m}{n}{<->s*[.
9]psyr}{} \DeclareSymbolFont{der@m}{U}{fsy}{m}{n}
\DeclareMathSymbol{\diff}{\mathord}{der@m}{182}

\newcommand{\qedd}{\hspace*{\fill}$\Box$\medskip}

\def\Y{{\mathbb{Y}}}
\def\U{{\mathbb{U}}}
\def\V{{\mathbb{V}}}
\def\L{{\mathbb{L}}}

\def\ff{{\mathcal F}}
\def\CI{{\mathcal{I}}}

\def\ee{{\mathcal E}}
\def\CJ{{\mathcal{J}}}

\def\P{{\mathbb P}}

\def\I{{\mathbb{I}}}    

\def\bu{{\mathbf{u}}}
\def\bv{{\mathbf{v}}}

\def\sat{\hbox{\rm{sat}}}
\def\asat{\hbox{\rm{asat}}}

\def\max{\hbox{\rm{max}}}

\def\deg{\hbox{\rm{deg}}}

\def\ord{\hbox{\rm{ord}}}
\def\initial{\hbox{\rm{I}}}

\def\sep{\hbox{\rm{S}}}

\def\dim{\hbox{\rm{dim}}}

\def\mod{\hbox{\rm{mod}}}
\def\leader{{\rm ld}}

\def\ord{\hbox{\rm{ord}}}

\def\h{\hbox{\rm{h}}}

\def\trdeg{\hbox{\rm{tr.deg}}}
\def\dtrdeg{\hbox{\rm{d.tr.deg}}}

\def\and{\cap}

\newcounter{bean}
\def\bl{\begin{list}{Step \arabic{bean}}{\usecounter{bean}}\labelwidth=34pt}
\def\el{\end{list}}

\def\deg{{\rm deg}}

\def\id{{\rm id}}

\def\normalization1{{\rm normalization1}}
\def\normalization{{\rm normalization}}

\def\irrfactor1{{\rm irrfactor1}}
\def\irrfactor{{\rm irrfactor}}

\def\sat{{\rm sat}}

\def\leader{{\rm ld}}

\begin{document}
\title{Differential Chow Form for Projective Differential Variety}
\author{Wei Li and Xiao-Shan Gao}
\address{KLMM, Academy of Mathematics and Systems Science\\
Chinese Academy of Sciences, Beijing 100190, China}
\email{xgao@mmrc.iss.ac.cn}
%
\thanks{Partially supported by a National Key Basic Research Project of China (2011CB302400) and
       by a Science Fund for Creative Research Groups from NSFC (60821002). }
\subjclass[2000]{Primary 12H05, 14C05; Secondary 14C17, 14Q99}

\date{July 16, 2011.}

\keywords{Differential Chow form; differentially homogenous prime
differential ideal; projective differential variety; intersection
theory.}

\begin{abstract}
In this paper, a generic intersection theorem in projective
differential algebraic geometry  is presented. Precisely, the
intersection of an irreducible projective differential variety of
dimension $d>0$ and order $h$ with a generic projective differential
hyperplane is shown to be an irreducible projective differential
variety of dimension $d-1$ and order $h$. Based on the generic
intersection theorem, the Chow form for an irreducible projective
differential variety is defined and most of the properties of the
differential Chow form in affine differential case are established
for its projective differential counterpart.  Finally, we apply the
differential Chow form to a result of linear dependence over
projective varieties given by Kolchin.

\end{abstract}

\maketitle

\section{introduction}

Differential algebra or differential algebraic geometry founded by
Ritt and Kolchin aims to study algebraic differential equations in a
similar way that polynomial equations are studied in commutative
algebra or algebraic geometry \cite{ritt,kol73}. Therefore, the
basic concepts of differential algebra geometry are based on those
of  algebraic geometry. An excellent survey on this subject can be
found in \cite{bc1}.

The Chow form, also known as the Cayley form or the Cayley-Chow
form, is a basic concept in algebraic geometry
\cite{gelfand,hodge,van} and has many important applications in
transcendental number theory  \cite{nes1,ph1}, elimination theory
\cite{brownawell,eisenbud,Pedersen,sturmfels}, algebraic
computational complexity \cite{Complexitychowform}.

Recently, the theory of differential Chow form for affine
differential algebraic varieties was developed \cite{gao}. It is
shown that most of the basic properties of algebraic Chow form can
be extended to its differential counterparts \cite{gao}.
Furthermore, a theory of differential resultant and sparse
differential resultant were also given  \cite{gao,sdresultant}.
In this paper, we will study differential Chow form for differential
projective varieties.

It is known that most results in projective algebraic geometry are
more complete than their affine analogs. But in the differential
case, this is not valid. Due to the complicated structure,
projective differential varieties are not studied thoroughly.
The basis of projective differential algebraic geometry was
established by Kolchin in his paper \cite{kol74}. There, he cited a
remark by Ritt:
\begin{quote}
 Consider an irreducible algebraic variety $V$ in complex
projective space $P_n(\mathbb{C})$ and $n+1$ meromorphic functions
$f_0,f_1,\ldots,f_n$ on some region of $\mathbb{C}$. J.F. Ritt once
remarked to me that there exists an irreducible ordinary
differential polynomial $h\in\mathbb{C}\{y_0,\ldots,y_n\}$,
dependent only on $V$ and having order equal to the dimension of
$V$, that enjoys the following property: A necessary and sufficient
condition that there exist $c_0,\ldots,c_n\in\mathbb{C}$ not all
zero such that $(c_0:\cdots:c_n)$ is a point of $V$ and
$\sum_{j}c_jf_j=0$ is that $(f_0,\ldots,f_n)$ be in the general
solution of $h$.
\end{quote}

Kolchin commented that ``This provides an occasion to describe the
beginning of a theory of algebraic differential equations in a
projective space $P_n(\ee)$."  And he devoted two papers on
differential projective spaces: \cite{kol74,kol92}, published
posthumously. In the former, Kolchin developed the foundation for a
theory of differentially homogenous differential ideals and their
differential zero sets in $P_n(\ee)$, which defines the Kolchin
topology in projective differential algebraic geometry. Also,
Kolchin proved Ritt's result mentioned above. In the following, we
will use ``the result on linear dependence over projective
varieties'' to refer to this result.

In this paper, we will first consider the dimension and order for
the intersection of a projective differential variety by a generic
projective differential hyperplane. Then we establish the theory of
differential Chow forms for projective differential varieties. As an
application, we will show that the differential polynomial $h$ in
Ritt's remark mentioned above  is just the differential Chow form of
the projective differential variety corresponding to $V$.

The rest of the paper is organized as follows. In section 2, we will
present the  basic notations and preliminary results in projective
differential algebraic geometry given by Kolchin which will  be used
in the paper. In section 3, we will prove the generic intersection
theorem of a projective differential variety by a generic projective
differential hyperplane. In section 4, the differential Chow form
for an irreducible projective differential variety is defined and
its basic properties are given. In section 5, we will apply the
differential Chow form theory to a result  on linear dependence over
projective varieties given by Kolchin. Finally, we present the
conclusion and propose a conjecture on the Chow form for the
projective differential variety.

\section{Preliminaries} In this paper, we fix $\ff$ to be an
ordinary differential field with derivation $\delta$, field of
constants $\mathcal {C}$, and universal differential field $\ee$
with field of constants $\mathcal {K}$. Primes and exponents $(i)$
indicate derivatives, and exponents $[t]$ indicate the set of
derivatives up to order $t$.

Let $\Y=(y_0,\ldots,y_n)$ and consider the differential polynomial
ring $\ff\{\Y\}=\ff\{y_0,\ldots,y_n\}$ over $\ff$. For any element
$t$ of any differential overring of $\ff\{\Y\}$, set
$t\Y=(ty_0,ty_1,\ldots,ty_n)$. Following Kolchin \cite{kol74}, we
introduce the concepts of differentially homogenous polynomials,
differentially homogenous differential ideals and differentially
projective varieties.
\begin{definition} \label{d-homogenous}
A differential polynomial $f \in \mathcal
{F}\{y_{0},y_{1},\ldots,y_{n}\}$ is called differentially homogenous
of degree $m$ if for a new differential indeterminate $t$, we have
$f(t y_{0},ty_{1}\ldots,t y_{n})=t^{m}f(y_{0},y_{1},\ldots,y_{n}) $.
\end{definition}

\begin{lemma} \cite{kol92} \label{le-separanthom}
Let $f$ be a nonzero differentially homogenous polynomial of degree
$d>0$ and let $o$ denote the order of $f$. Then for each index $j$,
$\partial f/\partial y_j^{(o)}$ is differentially homogenous of
degree $d-1$.
\end{lemma}

Now we will show that the initial and  the separant  of each
differentially homogenous polynomial w.r.t. any ranking is
differentially homogenous too.

\begin{theorem}\label{th-sepinihomo}
Let $f$ be a nonzero differentially homogenous polynomial of degree
$m$ and $\mathscr{R}$ be any ranking of $\Y$. Then both of its
initial and separant w.r.t. $\mathscr{R}$ are differentially
homogenous.
\end{theorem}

\proof We first consider the separant of $f$. Let $y_i^{(o)}$ be the
leader of $f$ w.r.t. $\mathscr{R}$. Then we have $t\frac{\partial
f}{\partial
y_i^{(o)}}(ty_0,\ldots,ty_n)=\frac{(ty_i)^{(o)}}{y_i^{(o)}}\frac{\partial
f}{\partial y_i^{(o)}}(ty_0,\ldots,ty_n)=\frac{\partial }{\partial
y_i^{(o)}}f(ty_0,\ldots,ty_n)=\frac{\partial}{\partial
y_i^{(o)}}t^mf(y_0,\ldots,y_n)=t^m\frac{\partial f}{\partial
y_i^{(o)}}(y_0,\ldots,y_n)$. Thus, $\frac{\partial f}{\partial
y_i^{(o)}}(t\Y)=t^{m-1}\frac{\partial f}{\partial y_i^{(o)}}(\Y)$.
It follows that the separant of $f$ is differentially homogenous of
degree $m-1$.

 Now we consider the initial of $f$.  Rewriting $f$ as a polynomial in $y_i^{(o)}$, we have
$f=I_f\cdot (y_i^{(o)})^l+I_1\cdot
(y_i^{(o)})^{l-1}+\cdots+I_{l-1}\cdot (y_i^{(o)})+I_l.$ As above, we
have proved that $\partial f/\partial y_i^{(o)}=lI_f\cdot
(y_i^{(o)})^{l-1}+(l-1)I_1\cdot (y_i^{(o)})^{l-2}+\cdots+I_{l-1}$ is
differentially homogenous with leader $y_i^{(o)}$. Continuing in
this way, $\frac{\partial^l f}{\partial (y_i^{(o)})^l}=l!I_f$ is
differentially homogenous. And it follows that $I_f$ is
differentially homogenous.\qedd

\hskip0.5cm

  More generally, let $(y_{ij})_{1\leq i\leq p,0\leq j\leq
n_i}$ be a family of differential indeterminates, set
$\Y_i=(y_{i0},\ldots,y_{in_i})\,(1\leq i\leq p)$, and consider the
differential polynomial ring
$\ff\{\Y_1,\ldots,\Y_p\}=\ff\{(y_{ij})_{1\leq i\leq p,0\leq j\leq
n_i}\}$ over $\ff$. Let $f\in\ff\{\Y_1,\ldots,\Y_p\}$ and
$(d_1,\ldots,d_p)\in \mathbb{N}^p$. If for any index $i$, $f$ is
differentially homogenous in $\Y_i$ of degree $d_i$, $f$ is said to
be {\em differentially $p$-homogenous} in $(\Y_1,\ldots,\Y_p)$ of
degree $(d_1,\ldots,d_p).$ When $\CI$ is a differential ideal of
$\ff\{\Y_1,\ldots,\Y_p\}$, then denote
$\CI:(\Y_1\cdots\Y_p)^\infty=\{f\in\ff\{\Y_1,\ldots,\Y_p\}\big|(y_{1j_1}\cdots
y_{pj_p})^ef\in\CI, 0\leq j_1\leq n_1,\ldots,0\leq j_p\leq
n_p\,\text{and for some $e$}\}$ and
$\CI:(\Y_1\cdots\Y_p)=\{f\in\ff\{\Y_1,\ldots,\Y_p\}\big|(y_{1j_1}\cdots
y_{pj_p})f\in\CI, 0\leq j_1\leq n_1,\ldots,0\leq j_p\leq n_p\}$.
Clearly,  $\CI:(\Y_1\cdots\Y_p)^\infty$ is a differential ideal, and
is a perfect ideal coinciding with $\CI:(\Y_1\cdots\Y_p)$ when $\CI$
is perfect.

\begin{definition}
Let $\CI$ be a differential ideal of $\ff\{\Y_1,\ldots,\Y_p\}$.
$\CI$ is called a differentially $p$-homogenous differential ideal
of $\ff\{\Y_1,\ldots,\Y_p\}$ if $\CI:(\Y_1\cdots\Y_p)=\CI$ and for
every $P\in\CI$ and a differential indeterminate $t$ over
$\ff\{\Y_1,\ldots,\Y_p\}$, $P(t\Y_1,\ldots,t\Y_p)\in\ff\{t\}\CI$ in
the differential ring $\ff\{t,\Y_1\cdots\Y_p\}$. In the case $p=1$
and $\Y_1=\Y$, $\CI$ is called a {\em differentially homogeneous
differential ideal} of $\ff\{\Y\}$.
\end{definition}

For $\CI\subset\ff\{\Y_1,\ldots,\Y_p\}$ and   any
 ranking $\mathscr{R}$ of $((y_{ij})_{1\leq i\leq
p,0\leq j\leq n_i}).$ $\CI$ has a characteristic set w.r.t.
$\mathscr{R}$ which is not unique. To impose the uniqueness
condition on characteristic set, Kolchin gave the definition of
canonical characteristic set which is unique for a differential
ideal and a fixed ranking $\mathscr{R}$. Let
$\Theta(\Y_1,\ldots,\Y_p)$ denote the set of all derivatives of
$y_{ij}^{(k)}\,(1\leq i\leq p,0\leq j\leq n_i,k\geq0)$ and let
$\textbf{M}$ denote the set of all differential monomials in
$\Y_1,\ldots,\Y_p$. Then the ranking $\mathscr{R}$ induces a lex
monomial order on $\textbf{M}$. For every nonzero $P\in
\ff\{\Y_1,\ldots,\Y_p\}$, $P$ is called unitary if the coefficient
of the highest differential monomial of $P$ is 1.
\begin{definition}
Let $\CI$ be a differential ideal of $\ff\{\Y_1,\ldots,\Y_p\}$ and
$\mathscr{R}$ be
 a ranking  of $(y_{ij}:1\leq i\leq p,0\leq
j\leq n_i)$. A characteristic set $\mathcal {A}$ of $\CI$ w.r.t.
$\mathscr{R}$ is said to be a {\em canonical characteristic set} if
it satisfies the following two conditions.

(i) For each element $A\in \mathcal {A}$ and every nonzero element
$B\in \CI$ that is reduced w.r.t. all the elements of $\mathcal {A}$
other than $A$, the highest differential monomial of $A$ is lower
than or equal to the highest differential monomial of $B$.

(ii) Each element of $\mathcal {A}$ is unitary.
\end{definition}
Kolchin proved that for a differential ideal and a fixed ranking
$\mathscr{R}$, $\CI$ has a unique canonical characteristic set.
Also, for a prime differential ideal of $\ff\{\Y_1,\ldots,\Y_p\}$,
according to the canonical characteristic set, Kolchin gave the
following theorem to test whether $\CI$ is differentially
$p$-homogenous.
\begin{theorem} \label{th-charset-homogenous}
Let $\CI$ be a prime differential ideal of $\ff\{\Y_1,\ldots,\Y_p\}$
and let $\mathcal{A}$ denote the canonical characteristic set of
$\CI$ w.r.t. some ranking of $(y_{ij})_{1\leq i\leq p,0\leq j\leq
n_i}.$ Then the followings   are equivalent.
\begin{itemize}

\item $\CI$ is differentially $p$-homogenous.
\item$\CI:(\Y_1\cdots\Y_p)=\CI$ and for every zero
$(\eta_1,\ldots,\eta_p)$ of $\CI$ in $\ee^{n_1+1}\times \cdots\times
\ee^{n_p+1}$ and each $s\in\ee\backslash\{0\}$ and each index $i$,
$(\eta_1,\ldots,s\eta_i,\ldots,\eta_p)$ is a zero of $\CI$.
\item
$\CI:(\Y_1\cdots\Y_p)=\CI$ and each element of $\mathcal{A}$ is
differentially $p$-homogenous in $(\Y_1,\ldots,\Y_p)$.

\end{itemize}
\end{theorem}

Let $n\in \mathbb{N}$ and consider the projective space
$\textbf{P}(n)$ over $\ee$. 
 Any element
$(a_0,a_1,\ldots,$ $a_n)$ of $\ee^{n+1}$ different from the origin
is a representative of a unique point $\alpha$ of $\textbf{P}(n)$,
denoted by $(a_0:a_1:\cdots:a_n)$. Given $n_1,\ldots,n_p\in
\mathbb{N}$, we consider the $p$-projective space
$\textbf{P}(n_1,\ldots,n_p)=\textbf{P}(n_1)\times
\cdots\times\textbf{P}(n_p)$. For any point
$\alpha=(\alpha_1,\ldots,\alpha_p)$ of $\textbf{P}(n_1,\ldots,n_p)$,
if $\textbf{a}_i=(a_{i0},a_{i1},\ldots,a_{in_i})$ is a
representative of $\alpha_i\,(1\leq i\leq p)$, then the element
$(a_{ij})_{1\leq i\leq p,0\leq j\leq
n_i}=(\textbf{a}_1,\ldots,\textbf{a}_p)$ of
$\ee^{n_1+1}\times\cdots\times\ee^{n_p+1}$ is called a
representative of $\alpha.$

Consider a differential polynomial $P\in\ee\{\Y_1,\ldots,\Y_p\}$ and a point $\alpha\in \textbf{P}(n_1,\ldots,$ $n_p)$.
 Say that $P$ vanishes at $\alpha$, and that $\alpha$ is a zero of $P$, if $P$ vanishes at every representative of $\alpha$.
 For a subset $\mathscr{M}$ of $\textbf{P}(n_1,\ldots,n_p)$, let $\I_\ff(\mathscr{M})$ denote the set of differential polynomials
 in $\ff\{\Y_1,\ldots,\Y_p\}$ that vanishes on  $\mathscr{M}$
 and write $\I(\mathscr{M})=\I_\ff(\mathscr{M})$.
Let $\V(S)$ denote the set of points of $\textbf{P}(n_1,\ldots,n_p)$
that are zeros of the subset $S$ of $\ee\{\Y_1,\ldots,\Y_p\}$. And a
subset $V$ of $\textbf{P}(n_1,\ldots,n_p)$ is called a 
{\em projective differential $\ff$-variety}
 if there exists $S\subset\ff\{\Y_1,\ldots,\Y_p\}$ 
 such that $V=\V(S)$.

As in the affine differential case, we have a one-to-one
correspondence between projective differential varieties and perfect
differentially homogenous differential ideals.
\begin{theorem}\cite{kol74}
The mapping from the set of projective differential $\ff$-varieties
of $\textbf{P}(n_1,$ $\ldots,$ $n_p)$ into the set of differentially
$p$-homogenous perfect differential ideals of
$\ff\{\Y_1,\ldots,\Y_p\}$, given by the formula
$V\longrightarrow\I_\ff(V)$, and the mapping in the opposite
direction, given by the formula $\CI\longrightarrow\V(\CI)$, are
bijective and inverse to each other. And  a projective differential
$\ff$-variety $ V$ is $\ff$-irreducible if and only if $\I(V)$ is
prime.
\end{theorem}

 At last, we give the following theorem on a property of
 differential specialization which will be used in the following.

\begin{theorem}\cite[Theorem 2.16]{gao}\label{th-specil}
Let $\{u_{1}, \ldots,  u_{r}\}\subset \ee$ be a set of differential
indeterminates over $\ff$,    and $P_{i}(\U, \Y)$ $\in \mathcal
{F}\{\U, \Y\}$ $(i=1, \ldots, m)$  differential polynomials in
$\U=(u_{1},\ldots,u_{r})$ and $\Y=(y_{1}, \ldots, y_{n})$. Let
$\overline{\Y}=(\overline{y}_{1}, \overline{y}_{2}, \ldots,
\overline{y}_{n})$,  where $\overline{y}_{i}\in \ee$ are
differentially free from $\ff\langle \U \rangle$. If $P_{i}(\U,
\overline{\Y})$ $(i=1, \ldots, m)$ are differentially dependent over
$\mathcal {F}\langle \U \rangle$, then for any specialization $\U$
to $\overline{\U}$ in $\mathcal {F}$, $P_{i}(\overline{\U},
\overline{\Y}) \, (i=1, \ldots,  m)$ are differentially dependent
over $\ff$.
\end{theorem}

\section{Generic Intersection Theory for  Projective Differential varieties }
In this section, we will consider the order and dimension of the
intersection of a projective differential variety by a generic
projective differential hyperplane. Before doing this, we first give
a rigorous definition of dimension and order for differentially
homogenous differential ideals.
\subsection{Order and Dimension in projective Differential Algebraic Geometry}
In the whole paper, when talking about prime differential ideals, we
always imply that they are distinct form the unit differential
ideal. For a differentially homogenous differential ideal
$\CI\subset \ff\{\Y\}\stackrel{\triangle}{=}\ff\{y_0,\ldots,y_n\}$,
we define differential independent set modulo $\CI$, parametric set
and differential dimension polynomial similar to the affine case.
More precisely, a variable set $\U\subset\Y$ is said to be an
independent set modulo $\mathcal {I}$ if $\mathcal {I}\cap\mathcal
{F}\{\U\}=\{0\}$. And a {\em parametric set} of $\mathcal {I}$ is a
maximal differentially independent set modulo $\mathcal {I}$. And
for a differentially homogenous prime differential ideal, we have
the following lemma.

\begin{lemma}\label{le-diffhomo-para}
Let $\CI$ be a  differentially homogenous prime differential ideal
in $\ff\{\Y\}$. Then its parametric set is not empty.
\end{lemma}
\proof Suppose the contrary. That is, for each $y_i$,
$\CI\cap\ff\{y_i\}\neq [0]\,(i=0,\ldots,n)$. Let $\mathcal {A}_i$ be
the canonical characteristic set of $\CI$ w.r.t. the elimination
ranking $y_i\prec y_0\prec\cdots\prec y_n$. Then by
Theorem~\ref{th-charset-homogenous}, each element of $\mathcal
{A}_i$ is differentially homogenous. Since $\CI\cap\ff\{y_i\}\neq
[0]$, there exists $A_{i0}\in\mathcal {A}_i$ such that
$A_{i0}\in\ff\{y_i\}$. If $\ord(A_{i0})>0$, it is easy to see that
$A_{i0}$ can not be differentially homogenous. So
$A_{i0}\in\ff[y_i]$. Using the fact that $\CI$ is prime,
$A_{i0}=y_i$ follows. Thus, for each $i$, $y_i\in\CI$. It follows
that $1\in \CI:\Y$. By Theorem~\ref{th-charset-homogenous}, we have
$\CI=\CI:\Y$, so $\CI=\ff\{\Y\}$, which is a contradiction. \qedd

%

 In \cite{kol64},
Kolchin introduced differential dimension polynomials for prime
differential polynomial ideals. Following Kolchin's definition, we
have the following definition for differentially homogenous prime
differential ideals.

\begin{definition}
Let $\CI$ be a differentially homogenous prime differential ideal of
$\ff\{\Y\}$. Then there exists a unique numerical polynomial
$\omega_{\CI}(t)$ such that
$\omega_{\CI}(t)=\dim(\CI\bigcap\ff[(y_i^{(k)})_{0\leq i\leq n,0\leq
k\leq t}])$ for all sufficiently big $t\in \mathbb{N}.$
$\omega_{\CI}(t)$ is called the {\em differential dimension
polynomial} of $\CI.$
\end{definition}

 For a differentially homogenous prime differential ideal of
$\ff\{\Y\}$, we can use the differential dimension polynomial to
give the definition of its differential dimension and order.
\begin{lemma}
Let $\CI$ be a differentially homogenous prime differential ideal of
$\ff\{\Y\}$. Then $\omega_{\CI}(t)$ can be written in the form
$\omega_{\CI}(t)=(d+1)(t+1)+o$ for $d\geq 0$. We define $d$ to be
the dimension of $\CI$ and $o$ is called the order of $\CI$.
\end{lemma}
\proof  By Lemma~\ref{le-diffhomo-para}, the cardinality of a
parametric set of $\CI$ is not less than 1. Then by \cite[Theorem
2]{kol64}, it follows.\qedd

\hskip0.5cm


In the affine case, we can read the information of a differential
ideal, such as dimension and order, from the perspective of its
generic point. However, in projective case,  it is a bit more
complex to do this. Firstly, we give the definition of generic
points for a differentially homogenous prime differential ideal
following Kolchin's notation\,(\cite{kol74}).

Consider a point $\alpha=(\alpha_1,\ldots,\alpha_p)\in
\textbf{P}(n_1,\ldots,n_p)$. Choose a representative
$\textbf{a}_i=(a_{i0},\ldots,a_{in_{i}})\in\ee^{n_i+1}\,(1\leq i\leq
p),$ and for each $i$ choose one index $j_i$ such that
$a_{ij_i}\neq0$. For any subfield $K$ of $\ee$, the field extension
$K((a_{ij_i}^{-1}a_{ij})_{1\leq i\leq p,0\leq j\leq n_i})$ is
independent of the choice of the representative $(a_0,\ldots,a_n)$
and indices $j_i$. Denote the field extension by $K(\alpha)$. When
$K(\alpha)=K$, the point $\alpha$ is said to be {\em rational} over
$K$. The set of points of $\textbf{P}(n_1,\ldots,n_p)$ that are
rational over $K$ is denoted by $\textbf{P}_K(n_1,\ldots,n_p)$.

Consider again the point $\alpha=(\alpha_1,\ldots,\alpha_p)\in
\textbf{P}(n_1,\ldots,n_p)$. Denote the differential field
$\ff\langle \ff( \alpha)\rangle$ by $\ff\langle\alpha\rangle$. 
Clearly, the differential transcendence polynomial of
$(a_{ij_i}^{-1}a_{ij})_{1\leq i\leq p,0\leq j\leq n_i}$ is
independent of the choices made above, and may therefore  be called
the differential transcendence polynomial of $\alpha$ over $\ff$,
denoted by $\omega_{\alpha/\ff}$. It can be written in the form
$\omega_{\alpha/\ff}(t)=a_1(t+1)+a_0$ where $a_i\in\mathbb{Z}$. Then
$a_1=\dtrdeg\,\ff\langle
\alpha\rangle/\ff$. 

Consider a second point $\alpha'=(\alpha_1',\ldots,\alpha_p')\in
\textbf{P}(n_1,\ldots,n_p)$ and a representative
$(\textbf{a}_1',\ldots,\textbf{a}_p')$ of $\alpha'$; for each $i$
write $\textbf{a}_i=(a_{i0},\ldots,a_{in_{i}})\in\ee^{n_i+1}$ and
fix $j_i'$ such that $a'_{ij_i'}\neq 0$. If $\I_{\ff}(\alpha)\subset
\I_{\ff}(\alpha')$, then $a_{ij'_{i}}\neq0\,(1\leq i\leq p)$, that
is, the indices $j_i$ can be chosen equal to the indices $j_i'$, and
evidently $(a'^{-1}_{ij'_{i}}a'_{ij})_{1\leq i\leq p,0\leq j\leq
n_{i}}$ is a differential specialization of
$(a^{-1}_{ij'_{i}}a_{ij})_{1\leq i\leq p,0\leq j\leq n_{i}}$ over
$\ff$. Conversely, if there exist indices $j_1,\ldots,j_p$ such that
$a_{ij_i}\neq0$, $a'_{ij_i}\neq0\,(1\leq i\leq p)$ and
$(a'^{-1}_{ij_{i}}a'_{ij})_{1\leq i\leq p,0\leq j\leq n_{i}}$ is a
differential specialization of $(a^{-1}_{ij_{i}}a_{ij})_{1\leq i\leq
p,0\leq j\leq n_{i}}$ over $\ff$, then $\I_{\ff}(\alpha)\subset
\I_{\ff}(\alpha')$. Under these conditions call $\alpha'$ a {\em
differential specialization of $\alpha$ over $\ff$}.

Let $\CI$ be a differentially $p$-homogenous prime differential
ideal of $\ff\{\Y_1,\ldots,\Y_p\}$ and $V$ be the corresponding
 projective differential $\ff$-variety of
 $\textbf{P}(n_1,\ldots,n_p)$; thus $V$ is $\ff$-irreducible,
 $V=\V_{n_1,\ldots,n_p}(\CI)$, and $\CI=\I(\V)$.
 And  for a point $\alpha\in\textbf{P}(n_1,\ldots,n_p)$,
$\I_{\ff}(\alpha)=\CI$ if and only if The set of all differential
specializations of $\alpha$ over $\ff $ is $V$.

Call such a point $\alpha$ a {\em generic point of $\CI$ in
$\textbf{P}(n_1,\ldots,n_p)$} or a {\em generic point of $V$ over
$\ff$}. We call $\omega_{\alpha/\ff}(t)$ the {\em differential
dimension polynomial} of $V$ and be denoted by $\omega_V$ and
$\dtrdeg\,\ff\langle\alpha\rangle/\ff$ may be called the {\em
differential dimension } of $V$. Then $\omega_V$ and $\omega_{\CI}$
have the following relation.

\begin{theorem} \label{th-diffdimpol}
Let $\CI$ be a differentially homogenous prime differential ideal of
$\ff\{\Y\}$ and $V=\V(\CI)$. Then
$\omega_{\CI}(t)=(t+1)+\omega_V(t)$.
\end{theorem}

\proof Without loss of generality, assume that $V\nsubseteq
\V(y_0)$. Let $(1,\xi_1,\ldots,\xi_n)$ be a generic point of $V$.
Let $u\in\ee$ be a differential indeterminate over
$\ff\langle\xi_1,\ldots,\xi_n\rangle$.
Firstly, we claim that
$\omega_{\CI}(t)=\trdeg\,\ff(u^{[t]},(u\xi_1)^{[t]},\ldots,(u\xi_n)^{[t]}
)/\ff$ for all sufficiently big $t\in\mathbb{N}$. Denote $\CI^a$ to
be the affine differential ideal in $\ff\{y_0,\ldots,y_n\}$
consisting of all elements of $\CI$. By the definition of
$\omega_{\CI}(t)$, it only needs to show that
$(u,u\xi_1,\ldots,u\xi_n)$ is a generic point of $\CI^a$. Firstly,
each polynomial in $\CI^a$ vanishes at $(u,u\xi_1,\ldots,u\xi_n)$.
Suppose $f^a\in\ff\{y_0,\ldots,y_n\}$ with
$f^a(u,u\xi_1,\ldots,u\xi_n)=0$. Since $u$ is a differential
indeterminate over $\ff\langle\xi_1,\ldots,\xi_n\rangle$, we can
regard $f^a(u,u\xi_1,\ldots,u\xi_n)$ as a differential polynomial in
$u$, which is identically zero. Thus, for each $\lambda\in\ee$,
$f^a(\lambda,\lambda\xi_1,\ldots,\lambda\xi_n)=0$. That is, $f^a$
vanishes at every representative of $(1,\xi_1,\ldots,\xi_n)$ and
$f^a\in \CI$ follows. So $f^a\in\CI^a$, and we have shown that
$(u,u\xi_1,\ldots,u\xi_n)$ is a generic point of $\CI^a$.

Thus, we have
\begin{eqnarray}\omega_{\CI}(t)&=& \nonumber\trdeg\,\ff(u^{[t]},(u\xi_1)^{[t]},\ldots,(u\xi_n)^{[t]}
)/\ff \\
\nonumber&=&\trdeg\,\ff(u^{[t]},(\xi_1)^{[t]},\ldots,(\xi_n)^{[t]}
)/\ff\\
\nonumber&=&\trdeg\,\ff(u^{[t]})/\ff+\trdeg\,\ff(u^{[t]})((\xi_1)^{[t]},\ldots,(\xi_n)^{[t]}
)/\ff(u^{[t]})\\
\nonumber&=&(t+1)+\trdeg\,\ff((\xi_1)^{[t]},\ldots,(\xi_n)^{[t]}
)/\ff\\
\nonumber&=&(t+1)+\omega_V(t).\nonumber
\end{eqnarray}
 \qedd
\begin{remark}
By the above theorem, we know that the differentially homogenous
prime differential ideal $\CI$ has the same dimension and order as
its corresponding projective differential variety.
\end{remark}

For every nontrivial differentially homogenous prime differential
ideal $\CI,$ there exists at least one index $i$ such that
$\CI\cap\ff\{y_i\}=[0]$. Of course, in this case, $y_i\notin\CI$.
Since we can permute variables when necessary, it only needs to
consider the case $\CI\cap\ff\{y_0\}=[0]$. Denote the set of all
nontrivial differentially homogenous prime differential ideals $\CI$
in $\ff\{y_0,y_1,\ldots,y_n\}$ with $\CI\cap\ff\{y_0\}=[0]$ by
$\mathcal {S}$. Denote the set of all prime differential ideals in
$\ff\{y_1,\ldots,y_n\}$ by $\mathcal {T}$. Now we give a one-to-one
correspondence between $\mathcal {S}$ and $\mathcal {T}$.

Define the map  $\phi: \,\mathcal
{S}\subset\ff\{y_0,y_1,\ldots,y_n\}\longrightarrow\, \mathcal
{T}\subset\ff\{y_1,\ldots,y_n\}$ and $\psi:\,\mathcal
{T}\subset\ff\{y_1,\ldots,y_n\}\longrightarrow\,\mathcal
{S}\subset\ff\{y_0,y_1,\ldots,y_n\}$ as follows: For each
$\CI\subset \mathcal{S}$, suppose $(\eta_0,\ldots,\eta_n)$ is a
generic point of $\CI$. Clearly, $\eta_0\neq0$. Let $\phi(\CI)$ be
the prime differential ideal in $\ff\{y_1,\ldots,y_n\}$ with
$(\eta_1/\eta_0,\ldots,\eta_n/\eta_0)$ its generic point.
Conversely, for a prime differential ideal $\CJ^a$ in
$\ff\{y_1,\ldots,y_n\}$ with a generic point $(\xi_1,\ldots,\xi_n)$,
let $\psi(\CJ^a)$ be a differentially homogenous prime differential
ideal in $\ff\{y_0,\ldots,y_n\}$ with $(v,v\xi_1,\ldots,v\xi_n)$ as
a generic point, where $v\in\ee$ is a differential indeterminate
over $\ff\langle\xi_1,\ldots,\xi_n\rangle$. Clearly, $\phi\circ
\psi=\id_{\mathcal {T}}$ and $\psi\circ \phi=\id_{\mathcal {S}}$. By
Theorem~\ref{th-diffdimpol},
$\omega_{\CI}(t)=(t+1)+\omega_{\phi(\CI)}(t)$. That is, $\CI$ has
the same dimension and order as $\phi(\CI).$

As above, we give the definition of $\phi$ and $\psi$ in the
language of generic points. Now we will give an interpretation from
the perspective of characteristic sets.

\begin{lemma} \label{le-phichar}
Let $\mathcal{A}:=A_1,\ldots,A_{n-d}$ be a canonical characteristic
set of $\CI\in\mathcal{S}$ w.r.t. any elimination ranking
$\mathscr{R}$ with $y_0\prec y_i$ for each $i$. Denote
$B_i=A_i(1,y_1,\ldots,y_n)$. Then $\mathcal{B}:=B_1,\ldots,B_{n-d}$
is a characteristic set of $\phi(\CI)$ w.r.t. the elimination
ranking induced by $\mathscr{R}$.
\end{lemma}

\proof Since $\CI$ is a differentially homogeneous prime
differential ideal, each $A_i$ is differentially homogenous by
Theorem~\ref{th-charset-homogenous}. Denote the separant and initial
of $A_i$ by $\sep_i$ and $\initial_i$. By
Theorem~\ref{th-sepinihomo}, $\sep_i$ and $\initial_i$ are
differentially homogenous. It follows that
$\sep_i(1,y_1,\ldots,y_n)$ and $\initial_i(1,y_1,\ldots,y_n)$ are
not zero. Consequently, $\leader(B_i)=\leader(A_i)$, and the
separant and initial of $B_i$ are $\sep_i(1,y_1,\ldots,y_n)$ and
$\initial_i(1,y_1,\ldots,y_n)$ respectively. So $\mathcal{B}$ is an
autoreduced set w.r.t. the elimination ranking induced by
$\mathscr{R}$.

Let $(\eta_0,\eta_1,\ldots,\eta_n)$ be a generic point of $\CI$.
Then $(\eta_1/\eta_0,\ldots,\eta_n/\eta_0)$ is a generic point of
$\phi(\CI)$. Clearly,
$B_i(\eta_1/\eta_0,\ldots,\eta_n/\eta_0)=A_i(1,\eta_1/\eta_0,\ldots,\eta_n/\eta_0)=0.$
That is, $B_i\in\phi(\CI)$, so $\mathcal{B}$ is an autoreduced set
in $\phi(\CI)$. Let $f$ be any polynomial in $\phi(\CI)$ and $r$ be
the remainder of $f$ w.r.t. $\mathcal{B}$. Then $r\in \phi(\CI)$ and
$r(\eta_1/\eta_0,\ldots,\eta_n/\eta_0)=0$. Let
$R=y_0^lr(y_1/y_0,\ldots,y_n/y_0)\in\ff\{\Y\}$. Then
$R(\eta_0,\eta_1,\ldots,\eta_n)=0$ and $R$ is reduced w.r.t.
$\mathcal{A}$. Thus, $R\equiv0$ and $r\equiv0$ follows. As a
consequence, $\mathcal{B}$ reduces every differential polynomial in
$\phi(\CI)$ to zero.
So, $\mathcal{B}$ is a characteristic set of $\phi(\CI)$ w.r.t. the
elimination ranking induced by $\mathscr{R}$. And the theorem is
proved. \qedd

\begin{remark}
Since the set of irreducible   projective differential varieties and
the set of differentially homogenous prime differential ideals have
a one-to-one correspondence, $\phi$ also gives a one-to-one map
between the set of  irreducible projective differential varieties
not contained in $y_0=0$ and the set of irreducible affine varieties
with inversive map $\psi$.
\end{remark}
\subsection{Generic Intersection theorem}

 In affine differential algebraic geometry, we have given the following
intersection theorem.

\begin{theorem}\cite[Theorem 3.14]{gao} \label{th-affineintersec}
Let $\mathcal {I}$ be a prime differential ideal in $\mathcal
{F}\{y_1,\ldots,y_n\}$ with dimension $d>0$ and order $h$. Let
$\{u_{0}, u_{1}, \ldots, u_{n}\}\subset \ee $ be a set of
differential indeterminates. Then $\mathcal {I}_1=[\mathcal {I},
u_{0}+u_{1}y_{1}+\cdots+u_{n}y_{n}]$ is a prime differential ideal
in $\mathcal {F}\langle u_{0}, u_{1}, \ldots,
u_{n}\rangle\{y_1,\ldots,y_n\}$ with dimension $d-1$ and order $h$.
\end{theorem}

Now we try to give the projective version of the above theorem.
Before doing this, we give the following lemma which will be used in
the proof of the main theorem.

\begin{lemma} \label{le-CI0}
Let $\mathcal {I}$ be a differentially homogenous prime differential
ideal in $\mathcal {F}\{\Y\}$ with dimension $d>0$ and
$u_0y_0+u_1y_1+\cdots+u_ny_n=0$ be a generic hyperplane. Denote
$\bu_0=(u_0,\ldots,u_n)$. Then the differential ideal
$\CI_0=[\CI,u_0y_0+u_1y_1+\cdots+u_ny_n]:(\Y)^\infty\subset\ff\{\Y;\bu_0\}$
is a differentially 2-homogenous prime differential ideal and
$\CI_0\cap\ff\{\bu_0\}=[0].$
\end{lemma}
\proof Let $(\xi_0,\ldots,\xi_n)$ be a generic point of $\CI$ that
is free from $\ff\langle u_0,\ldots,u_n\rangle$. Without loss of
generality, suppose $\xi_0\neq0$. Denote
$\CJ=[\CI,u_0y_0+\cdots+u_ny_n]:(\Y\bu)^\infty\subset\ff\{\Y;\bu\}$.
We prove that $\CI_0$ is a differentially 2-homogenous prime
differential ideal in $\ff\{\Y;\bu_0\}$ by showing that $\CJ$ is a
differentially 2-homogenous prime differential ideal in
$\ff\{\Y,\bu\}$ and $\CI_0=\CJ$.

Firstly, we show that for any point $\textbf{a} \in
\textbf{P}(n)\times \textbf{P}(n)$, if $\CJ$ vanishes at
$\textbf{a}$, then $\CJ$ vanishes at every representative of
$\textbf{a}$. Now, suppose $\CJ$ vanishes at $\textbf{a}$. For any
differential polynomial $H\in\CJ$, there exists some
$e\in\mathbb{N}$ such that $(y_ju_{k})^eH\in
[\CI,u_0y_0+\cdots+u_ny_n]$ for any $0\leq j,k\leq n$. Since $\CI$
and $u_0y_0+\cdots+u_ny_n$ vanish at every representative of
$\textbf{a}$, $H$ vanishes at it. It follows that $\CJ$ vanishes at
every representative of $\textbf{a}$. In this way, we say
$\textbf{a}$ is a zero of $\CJ$. Since $\CJ_1\subset\CJ$, $\CJ_1$
has the same property.

Now let
$\zeta=(\xi_0,\ldots,\xi_n;-(u_1\xi_1+\cdots+u_n\xi_n)/\xi_0,u_1,\ldots,u_n)$.
We now show that $\zeta$ is a generic point of $\CJ$. For any
$f\in\CJ$, there exists $e\in\mathbb{N}$ such that $(y_ju_{k})^eH\in
[\CI,u_0y_0+\cdots+u_ny_n]$. Take $j=0$ and $k=1$. Since $\xi_0\neq
0$, it follows that $H$ vanishes at $\zeta.$ Conversely, suppose $G$
is any differential polynomial in $\ff\{\Y;\bu\}$ such that $G$
vanishes at $\zeta$. Let $ G_1$ be the differential remainder of $G$
w.r.t. $u_0y_0+\cdots+u_ny_n$ with $u_0$ as its leader, then we have
$y_0^{a_0}G\equiv G_1,\,\mod\,[u_0y_0+\cdots+u_ny_n]$ for some
$a_0\in\mathbb{N}$. Then $G_1$ is free of $u_0$ and its derivatives
with $G_1(\zeta)=0$. Regard $G_1$ as a polynomial in $\ff\langle
u_1,\ldots,u_n\rangle$, then it vanishes at $(\xi_0,\ldots,\xi_n)$.
Since $[\CI]\subset\ff\langle u_1,\ldots,u_n\rangle \{\Y\}$ is a
prime differential ideal with a generic point
$(\xi_0,\ldots,\xi_n)$, $G_1\in[\CI]\cap\ff\{\Y;u_1,\ldots,u_n\}$.
Thus, $y_0^{a_0}G\in\CJ$. And for any index $j_0$ such that
$\xi_{j_0}\neq0$, similarly in this way, we can show that there
exists $a_{j_0}\in \mathbb{N}$ such that
$y_{j_0}^{a_{j_0}}G\in[\CI,u_0y_0+\cdots+u_ny_n]$. And if
$\xi_{j_0}=0,$ then $y_{j_0}\in\CI$, and
$y_{j_0}G\in[\CI,\P_0,\ldots,\P_d]$. Thus, it follows that $G\in
\CJ$ and $\zeta$ is a generic point of $\CJ$. Similarly, we can show
that $\zeta$ is also a generic point of $\CI_0$.
 By Theorem~\ref{th-charset-homogenous}, $\CI_0=\CJ$ is a
differentially $2$-homogenous prime differential ideal.

Suppose $\dim(\CI)>0.$ Then there exist $i\in\{1,\ldots,n\}$ such
that $\CI\cap\ff\{y_0,y_i\}=[0]$. It follows that $\xi_i/\xi_0$ is
differentially independent over $\ff$. By Theorem~\ref{th-specil},
$u_0,\ldots,u_n$ are differentially independent modulo $\CI_0$.
\qedd

\hskip0.5cm

Now, we give the following generic intersection theorem.
\begin{theorem}
Let $\mathcal {I}$ be a differentially homogenous prime differential
ideal in $\mathcal {F}\{\Y\}$ with dimension $d>0$ and order $h$.
Let $\{u_{0}, u_{1}, \ldots, u_{n}\}\subset \ee $ be a set of
differential indeterminates. Then $\mathcal {I}_1=[\mathcal {I},
u_{0}y_0+u_{1}y_{1}+\cdots+u_{n}y_{n}]:\Y^\infty$ is a
differentially homogenous prime differential ideal in $\mathcal
{F}\langle u_{0}, u_{1}, \ldots, u_{n}\rangle\{\Y\}$ with dimension
$d-1$ and order $h$.
\end{theorem}

\proof By Theorem~\ref{th-affineintersec},
$[\phi(\CI),u_{0}+u_{1}y_{1}+\cdots+u_{n}y_{n}]$ is a prime
differential ideal of dimension $d-1$ and order $h$. Notice the fact
that if we can show
$\phi(\CI_1)=[\phi(\CI),u_{0}+u_{1}y_{1}+\cdots+u_{n}y_{n}]$, then
it follows that $\CI_1$ is a differentially homogenous prime
differential ideal in $\mathcal {F}\langle u_{0}, u_{1}, \ldots,
u_{n}\rangle\{\Y\}$ with dimension $d-1$ and order $h$. So it
remains to show that
$\phi(\CI_1)=[\phi(\CI),u_{0}+u_{1}y_{1}+\cdots+u_{n}y_{n}]$.

Firstly, we show that $\CI_1$ is a differentially homogenous prime
differential ideal. Let $\CI_0$ be the differential ideal
$[\CI,u_{0}y_0+u_{1}y_{1}+\cdots+u_{n}y_{n}]:\Y^\infty$ in
$\ff\{\Y;\bu_0\}$ where $\bu=(u_0,\ldots,u_n)$. By
Lemma~\ref{le-CI0}, $\CI_0$ is a differentially 2-homogenous prime
differential ideal in $\ff\{\Y;\bu_0\}$ and
$\CI_0\cap\ff\{\bu_0\}=[0].$  Now we show that $\CI_1=[\CI_0]\subset
\ff\langle u_0,u_1,\ldots,u_n\rangle\{\Y\}$ is a nontrivial prime
differential ideal. If $1\in\CI_1$, then we have
$\CI_0\cap\ff\{u_0,u_1,\ldots,u_n\}\neq[0]$, a contradiction. So
$\CI_1\neq[1]$. Suppose $f_1,f_2\in\ff\langle
u_0,u_1,\ldots,u_n\rangle\{\Y\}$ with $f_1f_2\in\CI_1$. Then there
exist $h_i(\bu)\,(i=1,2)$ such that $h_if_i\in \ff\{\Y;\bu\}$. So
$(h_1f_1)(h_2f_2)\in\CI_0$. Then $h_1f_1\in\CI_0$ or
$h_2f_2\in\CI_0$, and it follows that $f_1\in \CI_1$ or $f_2\in
\CI_1$. Thus, $\CI_1$ is a differentially homogenous prime
differential ideal. Moreover,
 $\CI_1$ and $\CI_0$ have the
following relations:

(i)  $\CI_1\cap\ff\{\Y;\bu\}=\CI_0.$

(ii) Any characteristic set of $\CI_1$ with elements in
$\ff\{\Y;\bu\}$ such that $\bu$ are contained in the parametric set
is a characteristic set of $\CI_0.$

(iii)  Any characteristic set of $\CI_0$ with $\bu$ contained in its
parametric set is a characteristic set of $\CI_1$.

Since $\CI_1\neq[1]$, without loss of generality, we suppose
$\CI_1\cap\ff\{y_0\}=[0]$. Let $\mathcal{A}\subset\ff\{\Y;\bu\}$ be
a characteristic set of $\CI_1$ w.r.t. the elimination ranking
$y_0\prec y_1\prec\cdots\prec y_n.$ Then $\mathcal{A}$ is also a
characteristic set of $\CI_0$ w.r.t. the elimination ranking
$u_0\prec \cdots \prec u_n \prec y_0\prec y_1\prec\cdots\prec y_n$.
Let $\mathcal{B}$ be the autoreduced set obtained from $\mathcal{A}$
by setting $y_0=1$ in each element of $\mathcal{A}$. By
Lemma~\ref{le-phichar}, $\mathcal{B}$ is a characteristic set of
$\phi(\CI_1)$. Both $\phi(\CI_1)$ and
$[\phi(\CI),u_{0}+u_{1}y_{1}+\cdots+u_{n}y_{n}]$ are prime
differential ideals. If we can show that $\mathcal{B}$ is a
characteristic set of
$[\phi(\CI),u_{0}+u_{1}y_{1}+\cdots+u_{n}y_{n}]$, then
$\phi(\CI_1)=[\phi(\CI),u_{0}+u_{1}y_{1}+\cdots+u_{n}y_{n}]$
follows.

We claim that $\mathcal{B}$ is a characteristic set of
$\CI^a_0=[\phi(\CI),u_{0}+u_{1}y_{1}+\cdots+u_{n}y_{n}]\subset\ff\{\Y;\bu\}$.
We already know that if $(\xi_0,\ldots,\xi_n)$ is a generic point of
$\CI$, then $(\xi_0,\ldots,\xi_n;$
$-(u_1\xi_1+\cdots+u_n\xi_n)/\xi_0,u_1,\ldots,u_n)$ is a generic
point of $\CI_0$ and $(\xi_1/\xi_0,\ldots,\xi_n/\xi_0;$
$-(u_1\xi_1/\xi_0+\cdots+u_n\xi_n/\xi_0),u_1,\ldots,u_n)$ is a
generic point of $\CI^a_0$. Recall that $\mathcal{B}$ has the same
leaders as $\mathcal{A}$. It is easy to see that $B_i\in \CI^a_0$.
For any differential polynomial $f\in \CI^a_0,$ let $r$ be the
remainder of $f$ w.r.t. $\mathcal{B}$. Then
$r(\xi_1/\xi_0,\ldots,\xi_n/\xi_0;-(u_1\xi_1/\xi_0+\cdots+u_n\xi_n/\xi_0),u_1,\ldots,u_n)=0.$
Let $r^{\h}=y_0^lr(y_1/y_0,\ldots,y_n/y_0;u_0,\ldots,u_n)$ where $l$
is the denomination  of $r$. Then
$r^{\h}(\xi_0,\ldots,\xi_n;-(u_1\xi_1+\cdots+u_n\xi_n)/\xi_0,u_1,\ldots,u_n)=0$.
So $r^{\h}\in \CI_0$. But $\mathcal{A}$ is also a characteristic set
of $\CI_0$, thus, $r^{\h}=0$ and $r=0$ follows. Thus, $\mathcal{B}$
is a characteristic set of $\CI^a_0$. It follows that $\mathcal{B}$
is a characteristic set of
$[\phi(\CI),u_{0}+u_{1}y_{1}+\cdots+u_{n}y_{n}]$. Thus,
$\phi(\CI_1)=[\phi(\CI),u_{0}+u_{1}y_{1}+\cdots+u_{n}y_{n}]$ and the
theorem is proved. \qedd

As a corollary,  we give the following generic intersection theorem
in terms of projective differential variety.

\begin{theorem}
The intersection of an irreducible  projective differential variety
of dimension $d>0$ and order $h$ over $\mathcal {F}$ with
 a generic differential hyperplane  $
u_{0}y_0+u_{1}y_{1}+\cdots+u_{n}y_{n}=0$ is an irreducible
projective differential variety of dimension $d-1$ and order $h$
over $\ff\langle u_0,\ldots,u_n\rangle$.

\end{theorem}

\section{Differential Chow Forms for Projective Differential Varieties}
Let $V$ be an irreducible projective differential  $\ff$-variety in
$\textbf{P}(n)$ of differential dimension $d$ and order $h$. Suppose
$V$ does not lie in the hyperplane $y_0=0$. Let
$\xi=(1,\xi_1,\ldots,\xi_n)$ be a generic point of $V$. Let
$\bu=\{u_{ij}:\,i=0,\ldots,d;j=1,\ldots,n\}$ be a set of
differential indeterminates in $\ee$ over
$\ff\langle\xi_1,\ldots,\xi_n\rangle$. Denote $\zeta_i=-\sum_{j=1}^n
u_{ij}\xi_j\,(i=0,\ldots,d)$. Since $\dim(V)=d$,
$\dtrdeg\,\ff\langle\xi_1,\ldots,\xi_n\rangle/\ff=d$. Then as in the
affine differential case, we can prove that $\dtrdeg\,\ff\langle
\bu\rangle\langle\zeta_0,\ldots,\zeta_d\rangle/\ff\langle\bu\rangle=d$
and if $d>0$, then any  $d$ elements of $\zeta_i$ are differentially
independent over $\ff$. If $d=0$, then $\zeta_0$ is differentially
algebraic over $\ff.$ Let $u_{00},\ldots,u_{d0}$ be differential
indeterminates over $\ff\langle\bu,\xi_1,\ldots,\xi_n\rangle$. Since
$\zeta_0,\ldots,\zeta_d$ are differentially dependent over
$\ff\langle\bu\rangle$, there exists a relation
\[f(\bu;\zeta_0,\ldots,\zeta_d)=0\] where $f$ is a differential
polynomial in $\ff\langle\bu\rangle\{u_{00},\ldots,u_{d0}\}$ with
minimal order. We choose $f$  to be an irreducible differential
polynomial in $\ff\{\bu;u_{00},\ldots,u_{d0}\}$. Denote
$\bu_i=(u_{i0},\ldots,u_{in})\,(i=0,\ldots,d)$.
\begin{definition}\label{def-chowform}
The above irreducible differential polynomial
$f(\bu;u_{00},\ldots,u_{d0}) \in \ff\{\bu;u_{00},$ $\ldots,u_{d0}\}$
is defined to be {\em the differential Chow form} of $V$, denoted by
$F(\bu_0,\ldots,\bu_d)$.
\end{definition}

Let $\CI$ be the differentially homogenous prime differential ideal
in $\ff\{\Y\}$ associated to $V$, where $\Y=(y_0,\ldots,y_n)$. Let
$\P_i=u_{i0}y_0+u_{i1}y_1+\cdots+u_{in}y_n\,(i=0,\ldots,d)$. Then we
have the following theorem.
\begin{theorem} \label{th-chowhomogenous}
The ideal
$\mathcal{J}=[\CI,\P_0,\ldots,\P_d]:(\Y\bu_0\cdots\bu_d)^\infty$ in
$\ff\{\Y,\bu_0,\ldots,\bu_d\}$ is a differentially
$(d+2)$-homogenous prime differential ideal with $\mathcal
{J}\cap\ff\{\bu_0,\ldots,\bu_d\}=\sat(F)$, where $F$ is
differentially $(d+1)$-homogenous.
\end{theorem}
\proof Firstly, we show that for any point $\textbf{a} \in
\textbf{P}(n)\times \textbf{P}(n)\times\cdots
\times\textbf{P}(n)=\big(\textbf{P}(n))^{d+2}$, if $\CJ$ vanishes at
$\textbf{a}$, then $\CJ$ vanishes at every representative of
$\textbf{a}$. Now, suppose $\CJ$ vanishes at $\textbf{a}$.
 For any differential polynomial $H\in\CJ$, there
exists some $e\in\mathbb{N}$ such that $(y_ju_{0j_0}u_{1j_1}\cdots
u_{dj_d})^eH\in [\CI,\P_0,\ldots,\P_d]$ for any $0\leq
j,j_0,\ldots,j_d\leq n$. Since $\CI$ and $\P_0,\ldots,\P_d$ vanishes
at every representative of $\textbf{a}$, $H$ vanishes at it. It
follows that $\CJ$ vanishes at every representative of $\textbf{a}$.
In this way, we say $\textbf{a}$ is a zero of $\CJ$.

 To prove $\mathcal{J}$ is a prime differential ideal, it suffices
to show that $\textbf{c}=(1,\xi_1,\ldots,\xi_n;$
$\zeta_0,u_{01},\ldots,u_{0n};\ldots;\zeta_d,u_{d1},\ldots,u_{dn})$
is a generic zero of $\CJ$. Firstly, it is easy to see that
$\textbf{c}$ is a zero of $\CJ$. Now, suppose that $G$ is any
differential polynomial in $\ff\{\Y,\bu_0,\ldots,\bu_d\}$ such that
$G(\textbf{c})=0$. Now, $\P_0,\ldots,\P_d$ form an auto-reduced set
w.r.t. any elimination ranking of $\cup_{i=0}^d \bu_i\cup\Y$ such
that $z\prec u_{00}\prec\cdots \prec u_{d0}$ for any
$z\in\Y\cup\bu$. Let $G_1$ be the differential remainder of $G$
w.r.t. $\P_0,\ldots,\P_d$. Then $G_1$ is free of
$u_{i0}\,(i=0,\ldots,d)$ and
$$y_0^aG\equiv\,G_1,\,\mod\,[\P_0,\ldots,\P_d]\quad
(a\in\mathbb{N}).$$ So $G_1$ vanishes at $\textbf{c}$. Since
$G_1\in\ff\{\Y,\bu\}$, now rewritten $G_1$ as a differential
polynomial in $\bu$ with coefficients in $\Y$, i.e.,
$G_1=\sum_\phi\phi(\bu)g_\phi$ where $\phi(\bu)$ are different
differential monomials in $\bu$ and $g_\phi$ are differential
polynomials in $\ff\{\Y\}.$ Thus,
$G_1(\textbf{c})=\sum_\phi\phi(\bu)g_\phi(\xi)=0$. Since $\bu$ are
differential indeterminates over $\ff\langle \xi\rangle$,
$g_\phi(\xi)=0$ for any $\phi$. So, $g_\phi\in\CI$ and
$y_0^aG\in[\CI,\P_0,\ldots,\P_d]$.
And for any index $j_0$ such that $\xi_{j_0}\neq0$, it is easy to
show that $\zeta_0,\ldots,\zeta_d$ and
$\bu\backslash\{(u_{ij_0})_{0\leq i\leq d}\}$ are differentially
independent over $\ff\langle \xi\rangle$. Similarly in this way, we
can show that there exists $a_{j_0}\in \mathbb{N}$ such that
$y_{j_0}^{a_{j_0}}G\in[\CI,\P_0,\ldots,\P_d]$. And if $\xi_{j_0}=0,$
then $y_{j_0}\in\CI$, and $y_{j_0}G\in[\CI,\P_0,\ldots,\P_d]$. Thus,
it follows that $G\in \CJ$.

By Theorem~\ref{th-charset-homogenous}, $\CJ$ is a differentially
$(d+2)$-homogenous prime differential ideal. Clearly, $\mathcal
{J}\cap\ff\{\bu_0,\ldots,\bu_d\}$ is a differentially
$(d+1)$-homogenous prime differential ideal with a generic zero
$(\zeta_0,u_{01},\ldots,u_{0n};\ldots;\zeta_d,u_{d1},\ldots,u_{dn})$.
Since $\bu,\zeta_1,\ldots,\zeta_{d}$ are differentially independent
over $\ff$, the canonical characteristic set of $\mathcal
{J}\cap\ff\{\bu_0,\ldots,\bu_d\}$ consists of only one differential
polynomial, which is differentially $(d+1)$-homogenous by
Theorem~\ref{th-charset-homogenous}. By the definition of
differential Chow form above, this polynomial differs from $F$ by
only one factor  in $\ff.$ It follows that $\mathcal
{J}\cap\ff\{\bu_0,\ldots,\bu_d\}=\sat(F)$, and $F$ is differentially
$(d+1)$-homogenous. \qedd

By the above theorem, we have the following corollary.

\begin{lemma}
Let $F(\bu_{0},\bu_{1},\ldots,\bu_{d})$ be the differential Chow
form of an irreducible differential projective variety $V$ and
$F^*(\bu_{0},$ $\bu_{1},\ldots,$ $\bu_{d})$ obtained from $F$ by
interchanging $\bu_\rho$ and $\bu_\tau$. Then $F^*$ and $F$ differ
at most by a sign. Furthermore, $\ord(F,u_{ij})$ $(i=0,\ldots,d;
j=0,1,\ldots,n)$ are the same for all $u_{ij}$ occurring in $F$. In
particular, $u_{i0}\,(i=0,\ldots,d)$ appear effectively in $F$. And
a necessary and sufficient condition for some $u_{ij}\,(j>0)$ not
occurring effectively in $F$ is that $y_{j} \in \mathbb{I}(V).$
\end{lemma}
\proof By the definition of differential Chow form, similarly to
\cite[Lemma 4.9]{gao}, the lemma can be proved. \qedd

 From Definition~\ref{def-chowform}, we know that $F$ is also the differential Chow form
 of $\phi(V)$.  Since $V$ is of dimension $d$ and order $h$, $\phi
 (V)$ is of dimension $d$ and order $h$ too. Thus, by \cite[Theorem
4.12]{gao}, we have the following theorem.

\begin{theorem} \label{th-order}
Let $V$ be an irreducible  projective differential variety of
dimension $d$ and order $h$, and $F(\bu_0,\ldots,\bu_d)$ the
differential Chow form of $V$. Then $\ord(F)=h$.
\end{theorem}

Similarly to the affine case, the differential Chow form has the
following Poisson-type product formula.
\begin{theorem}
Let  $F(\bu_{0},\bu_{1},\ldots,\bu_{d})=f(\bu;u_{00},\ldots,u_{d0})$
 be the differential Chow form of an irreducible  projective differential $\ff$-variety of
dimension $d$ and order $h$. Then, there exist $\xi_{\tau
1},\ldots,\xi_{\tau n}$ in a differential extension field
$(\ff_\tau,\delta_\tau)$ ($\tau=1,\ldots,g$)  of $(\mathcal
{F}\langle\tilde{\bu}\rangle,\delta)$  such that
\begin{eqnarray}\label{eq-fac10}
F(\bu_{0},\bu_{1},\ldots,\bu_{d})=A(\bu_{0},\bu_{1},\ldots,\bu_{d})\prod^g_{\tau=1}(u_{00}+\sum_{\rho=1}^n
u_{0\rho}\xi_{\tau \rho})^{(h)}
\end{eqnarray}
where  $A(\bu_{0},\bu_{1},\ldots,\bu_{d})$ is in $\mathcal
{F}\{\bu_0,\ldots,\bu_d\}$, $\tilde{\bu}=\cup_{i=0}^d\bu_i\backslash
u_{00}$ and $g=\deg(f,u_{00}^{(h)})$.
Note that equation \bref{eq-fac10} is formal and should be
understood in the following precise meaning:
$(u_{00}+\sum_{\rho=1}^n u_{0\rho}\xi_{\tau \rho})^{(h)}
\stackrel{\triangle}{=}
\delta^{h}u_{00}+\delta_\tau^h(\sum_{\rho=1}^n u_{0\rho}\xi_{\tau
\rho}).$
\end{theorem}

For an element $\eta=(\eta_0,\eta_1,\ldots,\eta_n)$, denote its
truncation up to order $k$ as
$\eta^{[k]}=(\eta_0,\eta_1,\ldots,\eta_n,\ldots,\eta_0^{(k)},\eta_1^{(k)},$
$\ldots,\eta_n^{(k)})$.

Now we introduce the following notations:
\begin{equation}\label{eq-apoly}
 \begin{array}{l}
^a\P^{(0)}_0 =\,^a\P_0:=u_{00}y_0+u_{01}y_1+\cdots+u_{0n}y_n \\
^a\P^{(1)}_0 =\,
^a\P'_0:=u_{00}'y_0+u_{00}y_0'+u_{01}'y_1+u_{01}y'_1+\cdots+u_{0n}'y_n+u_{0n}y'_n
\\ \cdots \\
 ^a\P^{(s)}_0:=\sum_{j=0}^n \sum_{k=0}^s {s \choose k}u_{0j}^{(k)}y_j^{(s-k)}\end{array}
\end{equation}
which are considered to be algebraic polynomials in $\mathcal
{F}(\bu_{0}^{[s]},\ldots,\bu_{n}^{[s]})[\Y^{[s]}]$, and
$u_{ij}^{(k)}, y_i^{(j)}$ are treated as algebraic indeterminates.
A point $\eta=(\eta_0,\eta_1,\ldots,\eta_n)$ is said to be lying on
$^a\P^{(k)}_0$  if regarded as an algebraic point, $\eta^{[k]}$
%
is a zero of $^a\P^{(k)}_0$. Then as to the above Poisson product
formula, the following theorem holds.

\begin{theorem}
The points $(1,\xi_{\tau1},\ldots,\xi_{\tau n})\,(\tau=1,\ldots,g)$
in \bref{eq-fac10} are  generic points of the projective
differential
 $\ff$-variety $V$. If $d>0$, they also satisfy the
equations
\[\sum_{\rho=0}^n u_{\sigma \rho}y_{
\rho}=0\,(\sigma=1,\ldots,d).\] Moreover, they  are the only
elements of $V$ which also lie on $\P_i(i=1,\ldots,d)$\footnote{If
$d=0$, $\P_i(i=1,\ldots,d)$ is empty.} as well as on
$^a\P_0^{(j)}(j=0,\ldots,h-1)$.

\end{theorem}

As to the relations between the differential Chow form and the
projective differential variety, we have the following theorem.

\begin{theorem} \label{th-sf}
Let $F(\bu_{0},\bu_{1},\ldots,\bu_{d})$ be the differential Chow
form of $V$ and $S_{F}=\frac{\partial F}{\partial u_{00}^{(h)}}$.
Suppose that $\bu_i$ are differentially specialized over $\ff$ to
sets $\bv_i\subset\mathcal{E}$ and $\overline{\P}_{i}$ are obtained
by substituting $\bu_i$ by $\bv_i$ in $\P_i\,(i=0,\ldots,d)$.
If\, $\overline{\P}_i=0\,(i=0,\ldots,d)$ meet $V$, then
$(\bv_{0},\ldots,\bv_{d})$ is in the general solution  of the
differential equation $F=0$. Furthermore, if
$F(\bv_{0},\ldots,\bv_{d})=0$ and $S_{F}(\bv_{0},\ldots,\bv_{d})\neq
0$, then the $d+1$ differential hyperplanes $\overline{\P}_{i}=0$
$\,(i=0,\ldots,d)$ meet $V$.
\end{theorem}
\proof Let $\CI$ be the differentially homogenous prime differential
ideal in $\ff\{\Y\}$ associated to $V$. If
$\overline{\P}_i=0\,(i=0,\ldots,d)$ meet $V$, there exists
$\textbf{a}=(a_0,\ldots,a_n)\in\textbf{P}(n)$ with $a_{i_0}\neq0$
such that $\overline{\P}_i$ and $\CI$ vanish at $\textbf{a}.$
 Since
$[\CI,\P_0,\P_1,\ldots,\P_d]:(\Y)^\infty\cap\ff\{\bu_0,\bu_1,\ldots,\bu_d\}=\sat(F)$,
for any differential polynomial $G\in\sat(F)$, there exists
$e\in\mathbb{N}$ such that $y_j^eG\in[\CI,\P_0,\P_1,\ldots,\P_d]$
for every $j=0,1,\ldots,n$. So
$(a_{i_0})^eG(\bv_{0},\ldots,\bv_{d})=0$, and
$G(\bv_{0},\ldots,\bv_{d})=0$ follows. Thus,
$(\bv_{0},\ldots,\bv_{d})$ is in the general solution  of $F=0$.

Conversely,  suppose $F(\bv_{0},\ldots,\bv_{d})=0$ and
$S_{F}(\bv_{0},\ldots,\bv_{d})\neq 0$.  Let $(1,\xi_1,\ldots,\xi_n)$
be a generic point of $V$. Denote
$\zeta_i=-\sum_{j=1}^nu_{ij}\xi_j\,(i=0,\ldots,d)$. By the
definition of differential Chow form,
$F(\bu;\zeta_0,\ldots,\zeta_d)=0.$ Differentiating
$F(\bu;\zeta_0,\ldots,\zeta_d)=0$ w.r.t. $ u_{0\rho}^{(h)}$, we have
\begin{equation}\label{eq-parderiv}\overline{\frac{\partial F}{\partial
u_{0\rho}^{(s)}}}-\xi_{\rho} \overline{S_F}=0,\end{equation} where
$\overline{\frac{\partial F}{\partial u_{0\rho}^{(s)}}}$ and
$\overline{S_F}$ are obtained by replacing $(u_{00},\ldots,u_{d0})$
with $(\zeta_{0},\ldots,\zeta_{d})$ in $ \frac{\partial F}{\partial
u_{0\rho}^{(s)}}$ and $S_F$ respectively. Since
$(1,\xi_1,\ldots,\xi_n;\zeta_0,u_{01},\ldots,u_{0n};\ldots;\zeta_d,u_{d1},\ldots,u_{dn})$
is a generic point of $[\CI,\P_0,\P_1,\ldots,\P_d]:(\Y)^\infty$, by
equation~(\ref{eq-parderiv}), $A_\rho=S_Fy_{\rho}-\frac{\partial
F}{\partial u_{0\rho}^{(s)}}y_0\in
[\CI,\P_0,\P_1,\ldots,\P_d]:(\Y)^\infty$ ($\rho=1,\ldots,n$). It is
easy to see that $F,A_1,\ldots,A_n$ is a characteristic set of
$[\CI,\P_0,\P_1,\ldots,\P_d]:(\Y)^\infty$ w.r.t. the elimination
ranking $u_{01}\prec \cdots\prec u_{dn}\prec u_{10}\prec\cdots\prec
u_{d0}\prec u_{00}\prec y_0\prec y_1\prec\cdots \prec y_n$. Thus,
for any $H\in\CI$ and each $\P_i$,  there exist $e$ and $e_i$ such
that $S_F^e H\in [F,A_1,\ldots,A_n]$ and $S_F^{e_{i}} \P_i\in
[F,A_1,\ldots,A_n]$. Let $\bar{y}_i=\frac{\partial F}{\partial
u_{0\rho}^{(s)}}(\bv_{0},\ldots,\bv_{d})/S_F(\bv_{0},\ldots,\bv_{d})$
for $i=1,\ldots,n$. Clearly,
$(1,\bar{y}_1,\ldots,\bar{y}_n;\bv_{0},\ldots,\bv_{d})$ is a common
zero of $F,A_1,\ldots,A_n$. Thus,
$H(1,\bar{y}_1,\ldots,\bar{y}_n)=0$ and
$\overline{\P}_i(1,\bar{y}_1,\ldots,\bar{y}_n)=0.$ That is to say,
$(1,\bar{y}_1,\ldots,\bar{y}_n)$ is a common point of $V$ and
$\overline{\P}_i$. \qedd

\section{Application}

Let $V$ be an algebraic variety in $\textbf{P}_\mathcal {K}(n)$ that
is defined and irreducible over $\mathcal {C}$. Call an element
$\textbf{v}=(v_0,\ldots,v_n)\in\ee^{n+1}$ {\em linearly dependent
over} $V$ if there exists a point $\gamma\in V$ such that for some
(and hence every) representative $(c_0,\ldots,c_n)$ of $\gamma$,
$\sum_{j=0}^nc_jv_j=0$. If some representative of a point $\eta\in
\textbf{P}(n)$ is linearly dependent over $V,$ then every
representative of $\eta$ is. In this case, the {\em point} $\eta$ is
linearly dependent over $V$. In this section, $V$ will be fixed to
be an algebraic variety in $\textbf{P}_\mathcal {K}(n)$ that is
defined and irreducible over $\mathcal {C}$.

In \cite{kol74}, Kolchin gave the following theorem.

\begin{theorem} \label{th-RV}
Let $\mathfrak{R}$ denote the set of points of $\textbf{P}(n)$ that
are linearly dependent over $V$. Then there exists a differential
polynomial $R_V\in \mathcal{C}\{\Y\}$, irreducible over
$\mathcal{C}$, such that an element belong to $\mathfrak{R}$ if and
only if it is in the general solution of the differential equation
$R_V=0$. $R_V$ is unique up to a nonzero factor in $\mathcal{C}$ and
is differentially homogenous.
\end{theorem}

In a footnote at the end of his paper \cite{kol74}, Kolchin mentions
that H. Morikawa pointed out to him that the differentially
homogenous differential polynomial $R_V$ is the algebraic Chow form
of $V$ computed at the signed minors of the matrix
$(y^{(i)}_j)_{0\leq i\leq d,0\leq j\leq n}.$ That is, if
$G\big((u_{ij})_{0\leq i\leq d,0\leq j\leq n}\big)$ is the Chow form
of $V$, then $R_V=G\big((y_{j}^{(i)})_{0\leq i\leq d,0\leq j\leq
n}\big)$


\vskip10pt Clearly, $V$ has a natural structure of projective
differential $\mathcal{C}$-variety, defined by its polynomial
equations, together with the differential equations
$y_iy_j'-y_jy_i'=0\,(0\leq i<j\leq n)$, denoted by $V^\delta$. As we
have given the definition of differential Chow form in the last
section, in this section, we will explore the relationship between
$R_V$ and the differential Chow form of $V^\delta.$

\begin{lemma} \label{le-Vdelta}
If $\dim(V)=d,$ then $V^\delta$ defined as above is an irreducible
 projective differential variety of differential dimension 0 and
order $d$.
\end{lemma}
\proof Let $\CI_0$ be the defining ideal of $V$ in the polynomial
algebra $\mathcal{C}[\Y]=\mathcal[y_0,\ldots,y_n]$. Let
$(c_0,\ldots,c_n)\in\mathcal{K}^{n+1}$ be a generic point of
$\CI_0$. Without loss of generality, suppose $c_0\neq0$. Consider
the differential ideal $\CJ=[\CI_0,(y_iy_j'-y_jy_i')_{0\leq i<j\leq
n}]:(\Y)^\infty$ of $\mathcal{C}\{\Y\}$. We claim that
\begin{enumerate}
\item $\CJ$ is a differentially homogenous  prime differential ideal.
\item $V^\delta=\V(\CJ).$
\end{enumerate}
Suppose $u\in\ee$ is a differential indeterminate over
$\mathcal{K}$. To prove $\CJ$ is a prime differential ideal, it
suffices to prove that $(uc_0,\ldots,uc_n)$ is a generic point of
$\CJ$. Firstly, it is easy to show that $(uc_0,\ldots,uc_n)$ is a
zero of $\CJ$. Now, let $G\in\mathcal{C}\{\Y\}$ be a differential
polynomial vanishing at $(uc_0,\ldots,uc_n)$. Let $\mathscr{R}$ be
any ranking of $\Y$ such that $y_0\prec y_j\,(j=1,\ldots,n)$. Then
$\mathcal {A}:=y_0y_j'-y_jy_0'\,(j=1,\ldots,n)$ is an auto-reduced
set w.r.t. $\mathscr{R}$. Suppose the remainder of $G$ w.r.t.
$\mathcal {A}$ is $G_1$. Then $G_1$ is a differential polynomial in
$\mathcal{C}\{y_0\}[y_1,\ldots,y_n]$ and there exists some $a\in
\mathbb{N}$ such that $$y_0^aG\equiv G_1, \mod \,[\mathcal{A}].$$
Since $G(uc_0,\ldots,uc_n)=0$, $G_1(uc_0,\ldots,uc_n)=0$. Rewrite
$G_1$ as an algebraic polynomial in the proper derivatives of $y_0$
with coefficients in $\mathcal{C}[y_0,y_1,\ldots,y_n]$, then we have
$G_1=\sum_\phi\phi(y_0',y_0'',\ldots)G_\phi(y_0,y_1,\ldots,y_n)$
where $\phi(y_0',y_0'',\ldots)$ are distinct  monomials in
$y_0',y_0'',\ldots$. Since $u',u'',\ldots$ are algebraic
indeterminates over $\mathcal{K}(u)$, $G_\phi(uc_0,\ldots,uc_n)=0,$
so, $G_\phi\in \CI_0$ for each $\phi$. Thus,
$y_0^aG\in[\CI_0,(y_iy_j'-y_jy_i')_{0\leq i<j\leq n}]$. Similarly in
this way, for any $j_0$ such that $c_{j_0}\neq0$, we can show that
$y_{j_0}^{a_{j_0}}G\in[\CI_0,(y_iy_j'-y_jy_i')_{0\leq i<j\leq n}]$
for some $a_{j_0}\in \mathbb{N}$. So, $G\in \CJ$ and
$(uc_0,\ldots,uc_n)$ is a generic point of $\CJ$. Thus, $\CJ$ is a
prime differential ideal. Clearly, $\CJ:(\Y)=\CJ$ and for any zero
$\eta$ of $\CJ$ and every $s\in \ee^\star$, $s\eta$ is a zero of
$\CJ$. By Theorem~\ref{th-charset-homogenous}, $\CJ$ is a
differentially homogenous prime differential ideal. And it is easy
to show that $V^\delta=\V(\CJ)$. So $V^\delta$ is an irreducible
differentially projective variety and $(uc_0,\ldots,uc_n)$ is a
generic point of it.

The differential dimension polynomial of $V^\delta$ is
\begin{eqnarray}\omega_{V^\delta}(t)&=&\trdeg\,\mathcal{C}\Big(\big(\frac{uc_j}{uc_0}\big)^{(k)}:1\leq j\leq
n,k\leq
t\Big)\big/\mathcal{C} \nonumber \\
&=&\trdeg\,\mathcal{C}\big((c_j/c_0)^{(k)}:1\leq j\leq n,k\leq
t\big)\big/\mathcal{C}\nonumber\\
&=&\trdeg\,\mathcal{C}\big(c_j/c_0:1\leq j\leq
n\big)\big/\mathcal{C}\nonumber\\ &=& \dim(V)=d.
\nonumber\end{eqnarray} Thus, $V^\delta$ is of dimension $0$ and
order $d$ and the lemma follows. \qedd

Now we give the main theorem as follows.

\begin{theorem}\label{th-r}
The differential polynomial $R_V(\bu_0)$ defined in
Theorem~\ref{th-RV} is equal to the differential Chow form of
$V^\delta$ in the sense of multiplied by a nonzero constant in
$\mathcal{C}$.
\end{theorem}
\proof By Lemma~\ref{le-Vdelta}, $V^\delta$ is an irreducible
projective differential variety of dimension 0 and order $d$. Let
$F(\bu_0)$ be the differential Chow form of $V^\delta$. Then by
Theorem~\ref{th-order}, $\ord(F,\bu_0)=d$. By Theorem~\ref{th-RV},
$\ord(R_V,\bu_0)=d$.

Since $R_V(\bu_0)$ and $F(\bu_0)$ are irreducible differential
polynomials in $\mathcal{C}\{\bu_0\}$ with the same order, if we can
prove
$R_V(\bu_0)\in[\I(V^\delta),\P_0]:(\Y\bu_0)^\infty=\sat(F(\bu_0))$,
then by Theorem~\ref{th-chowhomogenous}, $R_V$ and $F$ differs at
most by a nonzero constant in $\mathcal{C}.$ Now we are going to
show that $R_V(\bu_0)\in[\I(V^\delta),\P_0]:(\Y\bu_0)^\infty$. In
the proof, we shall follow the notations in the above lemma. For
instance, $\CI_0=\I(V)\subset\mathcal{C}\{\Y\}$ and
$(c_0,\ldots,c_n)$ is a generic point of $V.$

Let $F_0(\bu_0,\ldots,\bu_d)$ be the algebraic Chow form of $V$,
where $\bu_1,\ldots,\bu_d$ are the vectors of coefficients of the
generic algebraic hyperplanes $\L_1,\ldots,\L_d$ respectively, and
$\P_0$ is regarded as an algebraic hyperplane at the very moment. By
the definition of algebraic Chow form,
$(\CI_0,\P_0,\L_1,\ldots,\L_d):\Y^\infty=\asat(F_0).$ Suppose
$c_{j_0}\neq0$. Then there exists $a_{j_0}$ such that
\begin{equation}\label{eq-0}y_{j_0}^{a_{j_0}}F_0(\bu_0,\ldots,\bu_d)=\sum_{k=0}^d
h_k\L_k+\sum_{i}g_if_i\end{equation} where
$f_i\in\CI_0\subset\mathcal{C}[\Y]$ and $h_k,g_i
\in\mathcal{C}[\Y,\bu_0,\ldots,\bu_d].$ Denote
$\bu_0^{(k)}=(u_{00}^{(k)},\ldots,u_{0n}^{(k)})$. Replace $\bu_k$ by
$\bu_0^{(k)}$ in (\ref{eq-0}) for $k=1,\ldots,d$, we obtain
\begin{equation}\label{eq-1}y_{j_0}^{a_{j_0}}F_0(\bu_0,\bu_0',\ldots,\bu_0^{(d)})=\sum_{k=0}^d \widehat{h_k}\widehat{\L_k}+\sum_{i}\widehat{g_i}f_i,\end{equation}
where
$\widehat{\L_k}=u_{00}^{(k)}y_0+u_{01}^{(k)}y_1+\ldots+u_{0n}^{(k)}y_n$.
Denote $G_{ij}=y_iy_j'-y_jy_i'\,(0\leq i<j\leq n)$ and for $i>j,
\,G_{ij}=-G_{ji}$. Now we claim that there exists $b_k\in
\mathbb{N}$ such that
$$y_{j_0}^{b_k}\widehat{\L_k}\in[\P_0,(G_{ij})_{0\leq i<j\leq
n}],\qquad (k=0,\ldots,d).$$ Assuming the claim holds. Denote
$e_{j_0}=\max_{k}\{a_{j_0}+b_{k}\}$. Then by (\ref{eq-1}), for each
$j_0$ such that $c_{j_0}\neq0$, we have
$y_{j_0}^{e_{j_0}}F_0(\bu_0,\bu_0',\ldots,\bu_0^{(d)})\in
[\I(V^\delta),\P_0]$. So
$R_V(\bu_0)=F_0(\bu_0,\bu_0',\ldots,\bu_0^{(d)})\in
[\I(V^\delta),\P_0]:(\Y\bu_0)^\infty$ follows. Thus, it suffices to
prove the claim.

Before proving the claim, we first show that for each $y_l$ and each
$m$, $y_{j_0}^{m-1}y_l^{(m)}\equiv h_m(y_{j_0})\cdot y_l,\,\mod
\,\,[G_{j_0l}]$ where $h(y_{j_0})\in\mathcal{C}\{y_{j_0}\}$. For
$m=0$, it is trivial and $m=1$, $y_{j_0}y_l'=y_{j_0}'y_l+G_{j_0l}$.
Suppose it holds for $0,\ldots,m-1$. Since
$G_{j_0l}^{(m-1)}=\big(y_{j_0}y_l'-y_{j_0}'y_l\big)^{(m-1)}=\sum_{s=0}^{m-1}{m-1\choose
s}y_{j_0}^{(m-1-s)}y_l^{(s+1)}-\sum_{s=0}^{m-1}{m-1\choose
s}y_{j_0}^{(m-s)}y_l^{(s)}$,
$y_{j_0}y_l^{(m)}=G_{j_0l}^{(m-1)}+\sum_{s=1}^{m-1}\big({m-1\choose
s}-{m-1\choose s-1}\big)y_{j_0}^{(m-s)}y_l^{(s)}+y_{j_0}^{(m)}y_l$.
By the hypothesis, $y_{j_0}^{s-1}y_l^{(s)}\equiv
h_s(y_{j_0})y_l,\,\mod \,[G_{j_0l}]$ holds for $s\leq m-1$. Thus,
$y_{j_0}^{m-1}y_l^{(m)}=y_{j_0}^{m-2}G_{j_0l}^{(m-1)}+\sum_{s=1}^{m-1}\big({m-1\choose
s}-{m-1\choose
s-1}\big)y_{j_0}^{(m-s)}y_{j_0}^{m-2}y_l^{(s)}+y_{j_0}^{(m)}y_{j_0}^{m-2}y_l\equiv
h_m(y_{j_0})\cdot y_{l},\mod\,[G_{j_0l}]$, where
$h_m=\sum_{s=1}^{m-1}\big({m-1\choose s}-{m-1\choose
s-1}\big)y_{j_0}^{(m-s)}y_{j_0}^{m-s-1}h_s+
y_{j_0}^{(m)}y_{j_0}^{m-2}$.

Now we are going to prove the  claim by induction on $k$. Firstly,
$\L_0=\P_0$. And for $k=1$,
$\widehat{\L_1}=\sum_{l=0}^nu_{0l}'y_l=(\sum_{l=0}^nu_{0l}y_l)'-\sum_{l=0}^nu_{0l}y_l'=\P_0'-\sum_{l=0}^nu_{0l}y_l'.$
Then
$y_{j_0}\widehat{\L_1}=y_{j_0}\P_0'-\sum_{l=0}^nu_{0l}y_{j_0}y_l'=y_{j_0}\P_0'-\sum_{l=0}^nu_{0l}\big(G_{j_0l}-y_{j_0}'y_l\big)
=y_{j_0}\P_0'-\sum_{l=0}^nu_{0l}G_{j_0l}+y_{j_0}'\P_0$. So it holds
for $k=1$. Now suppose the claim holds for integers less than $k$
and we shall deal with $k$. Since
$\widehat{\L_k}=\sum_{l=0}^nu_{0l}^{(k)}y_l=\P_0^{(k)}-\sum_{l=0}^n\sum_{m=1}^{k}{k\choose
m}u_{0l}^{(k-m)}y_l^{(m)}$,
\begin{eqnarray}
y_{j_0}^{k-1}\widehat{\L_k}&=&y_{j_0}^{k-1}\P_0^{(k)}-\sum_{m=1}^{k}{k\choose
m}\sum_{l=0}^n u_{0l}^{(k-m)}y_{j_0}^{k-1}y_l^{(m)} \nonumber\\
&=& y_{j_0}^{k-1}\P_0^{(k)}-\sum_{m=1}^{k}{k\choose
m}h_m(y_{j_0})y_{j_0}^{k-m}\big(\sum_{l=0}^n
u_{0l}^{(k-m)}y_l\big),\,\mod\, [(G_{ij})_{0\leq i<j\leq n}].
\nonumber
\end{eqnarray}
 By
hypothesis, for $s<k$, $y_{j_0}^{b_s}\widehat{\L_s}\in
[\P_0,(G_{ij})_{0\leq i<j\leq n}]$. Thus, the claim follows.

 \qedd

Combining Theorem~\ref{th-RV} with Theorem~\ref{th-r}, we have the
following corollary.

\begin{cor}\label{cor-chownecesuff}
Let $V^\delta$ be defined as above and
$F(u_{00},u_{01},\ldots,u_{0n})$ be its projective differential Chow
form. If $v_{0j}\in\mathcal {E}$ is any specialization of
$u_{0j}\,(j=0,1,\ldots,n)$, then $V^\delta$ and the differential
hyperplane
$$u_{00}y_0+u_{01}y_1+\cdots+u_{0n}y_n=0$$ have points in common if and only if  $(v_{00},v_{01},\ldots,v_{0n})$
 is in the general solution of the differential equation $F=0$.
\end{cor}

\section{Conclusion}

In this paper, we first prove a theorem for the intersection of an
irreducible projective differential variety with  generic projective
differential hyperplanes. Then we define the  differential Chow form
for projective differential varieties and give its basic properties.
Finally, we show that the formula  on linear dependence over an
algebraic projective variety given by Kolchin is actually the
differential Chow form of the projective variety teated as a
differential projective variety in certain sense.

For an algebraic projective variety $V$ of dimension $d$, its
projective Chow form gives a necessary and sufficient condition for
$V$ having common points with the $d$ hyperplanes
$\sum_{j=0}^nu_{ij}x_j=0$ \cite[p.50]{hodge}. For a particular kind
of projective differential varieties, Corollary
\ref{cor-chownecesuff} tells us that $\sat(F)$ gives a necessary and
sufficient condition that $V$ and the differential hyperplane
$\sum_{j=0}^nu_{0j}y_j=0$ intersect, where $F$ is the projective
differential Chow form of $V$. However, up to now, for general
projective differential varieties, Theorem~\ref{th-sf} only gives a
necessary condition. We conjecture here that for general projective
varieties, the differential saturation ideal of its differential
Chow form also gives such a necessary and sufficient condition.

\end{document}